\documentclass[12pt,twoside]{article}
\usepackage{amsmath, amsfonts, amssymb, amsthm,latexsym}
\usepackage{amscd,graphics}
\usepackage[cp1251]{inputenc}
\usepackage[english]{babel}

\newtheorem{theorem}{Theorem}[section]
\newtheorem{lemma}[theorem]{Lemma}
\newtheorem{corollary}[theorem]{Corollary}
\newtheorem{proposition}[theorem]{Proposition}

\newtheorem{remark}{Remark}

\newtheorem{hypothesis}{Hypothesis}

\begin{document}

\newcommand{\iiRGH}{{\rm Rep}\,(G,{\mathcal H})}
\newcommand{\iiRGX}{{\rm Rep}\,(G,X)}
\newcommand{\iiRGdf}{{\rm Rep}\,(G,d,f)}
\newcommand{\iiRGU}{{\rm Rep}\,(G,\coprod')}
\newcommand{\iiRGUC}{{\rm Rep}_\circ\,(G,\coprod')}
\newcommand{\iiRGUB}{{\rm Rep}_\bullet\,(G,\coprod')}
\newcommand{\iiRGUn}[1]{{\rm Rep}^{#1}\,(G,\coprod')}
\newcommand{\iiRGdfC}{{\rm Rep}_\circ\,(G,d,f)}
\newcommand{\iiRGdfB}{{\rm Rep}_\bullet\,(G,d,f)}
\newcommand{\iiSC}[1]{\stackrel{\circ}{#1}}
\newcommand{\iiSB}[1]{\stackrel{\bullet}{#1}}
\newcommand{\iiSCB}[1]{\stackrel{\circ\bullet}{#1}}
\newcommand{\iiSBC}[1]{\stackrel{\bullet\circ}{#1}}
\newcommand{\iiCor}[2]{\langle #1 \rangle^{#2}}

%\noindent UDK 519.1

%\setcounter{page}{20}

\begin{center}
\LARGE \bf
Singular Locally-Scalar \\
Representation of Quivers \\
in Hilbert Spaces \\
and Separating Functions

\end{center}

\begin{center}
\Large\bfseries I.K. Redchuk~$*$, A.V. Roiter~$^{**}$
\end{center}

{ \footnotesize

$^{*}$~Institute of Mathematics of National Academy of Sciences of Ukraine,\\
Tereshchenkovska str., 3, Kiev, Ukraine, ind. 01601\\
E-mail: red@imath.kiev.ua
\medskip

$^{**}$~Institute of Mathematics of National Academy of Sciences of Ukraine,\\
Tereshchenkovska str., 3, Kiev, Ukraine, ind. 01601\\
E-mail: roiter@imath.kiev.ua}

\bigskip

In the classical works \cite{Gabr}, \cite{BerGelPon}
representations of quivers (in the category of finite-dimensional vector
spaces) were considered and their connection with root systems
of corresponding graphs was established, which developing further
in the works \cite{Naz2}--\cite{GabRoi}. Repeatedly attempts were
made to generalize representations of quivers to metric and, in particular,
Hilbert spaces, but at that mentioned connection was lost. In
\cite{KrugRoit} a restriction of {\it local scalarity} was
imposed on representations of graphs in Hilbert spaces, and after that
it was managed to generalize the results of
\cite{Gabr}, \cite{BerGelPon} in a natural way by constructing, analogous
to Coxeter functors, functors of even and add reflections, moreover
it turned out that just these representations are of interest of
functional analysis and in some particular cases in other terms
in fact were considered in \cite{OstSam}--\cite{KruRab}.

In this paper we will consider a connection of locally-scalar
representations of extended Dynkin graphs with function $\rho$
considered in \cite{NazRoit}, and also functions $\rho_k$ playing analogous
part for more wide class of graphs.

\section{Standard singular representations of extended Dynkin graphs.}

We will widely use denotations and definitions of the article
\cite{KrugRoit}.

All considered graphs will be supposed finite, connected and
acyclic (i.~e. woods).

{\it Multiplicity} $\mu(g)$ of a vertex $g \in G_v$ is $|M_g|$ ($M_g =
\{g_i \in G_v\, |\, g_i-g\}$ \cite{KrugRoit}).

\begin{hypothesis}\label{h1}

All indecomposable locally-scalar representation of a graph $G$ are
finite-dimensional, if and only if $G$ is a Dynkin graph or extended
Dynkin graph.

\end{hypothesis}

In this paper we will consider only finite-dimensional
representations.

Let us fix a separation of set of the vertices $G_v$ of graph $G$ to even
$\iiSC{G_v}$ and odd ones $\iiSB{G_v}$ \cite{KrugRoit} and a numeration
$g_1,\ldots ,g_n$, numbering (in an arbitrary order) at first odd
$g_1,g_2, \ldots ,g_p$ and then even vertices $g_{p+1},g_{p+2}, \ldots
,g_n$. Let $x \in V_G$, ($x : G_v \to \mathbb{C}$) $x_i = x(g_i)$,
$c$~--- Coxeter transformation on $V_G$, $c = \sigma_{g_n}
\sigma_{g_{n-1}} \cdots \sigma_{g_1}$, $(\sigma_{g_i}(x))_i = - x_i +
\sum\limits_{j,\; g_j \in M_{g_i}} x_j$, $(\sigma_{g_i}(x))_j = x_j$ when $j
\neq i$. It is clear that $\sigma_i^2 = \rm id$ $(i = \overline{1,n})$.
Therefore $c^{-1} = \sigma_{g_1} \cdots \sigma_{g_n}$.

Vector $x \in V_G^+$ is {\it regular} if $c^t(x) \in V_G^+$
for any $t \in \mathbb{Z}$ and {\it singular} in a contrary case
(terminology traces back to~\cite{GelPon}).

The main part in our study will play extended Dynkin graphs
($\widetilde{D}_n$, $\widetilde{E}_6$, $\widetilde{E}_7$,
$\widetilde{E}_8$). Their connection with Coxeter transformations can
be characterized by the following well-known statement.

\begin{lemma}\label{l_u_EGD}

If $G$ is a wood, $u_G = u \in V_G^+$ and $c(u) = u$, then $G$ is an
extended Dynkin graph and (up to the common multiplier) $u_G$ looks like:

\begin{center}

Graph $\widetilde{D}_n$

\includegraphics{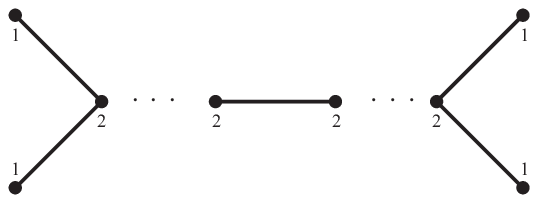}

\medskip

Graph $\widetilde{E}_6$

\includegraphics{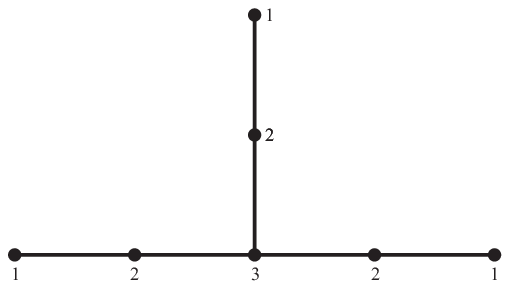}

\medskip

Graph $\widetilde{E}_7$

\includegraphics{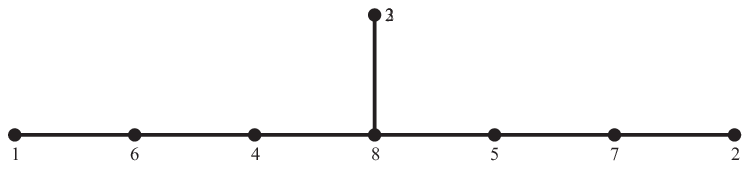}

\medskip

Graph $\widetilde{E}_8$

\includegraphics{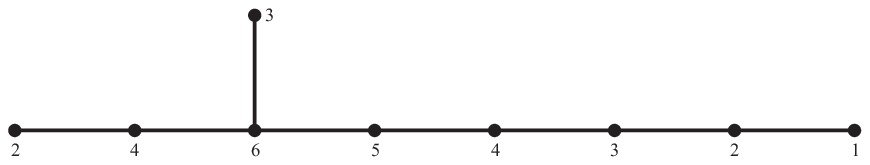}

\end{center}

\end{lemma}

{\bf Proof.}

$c(u) = u$ implies $\sigma_{g_i}(u) = u$ (and this is equivalent to the
statement, that when $i = \overline{1,n}$ $2u(g_i) = \sum\limits_{g \in
M_{g_i}}u(g)$) and $u_i > 0$ when $i \in \overline{1,n}$. Vertex $g$ is
the {\it point of branching} if $\mu(g) > 2$; $g$~--- {\it point of weak
branching} if $\mu(g) = 3$ and $|M_g| = \{a,b,c\}$ where $\mu(a) = \mu(b) =
1$. $\sigma_{g_i} = u$ implies:

1) if $a-b$ then $u(a) \geq \frac{1}{2}u(b)$, at that $u(a) =
\frac{1}{2}u(b)$ only if $\mu(a) = 1$;

2) if $\mu(g) > 2$, $g-a$ then $u(a) \leq u(g)$, at that $u(g) = u(a)$
only if $M_g = \{a,b,c\}$ and $\mu(b) = \mu(c) = 1$ (it follows from 1));

3) if $a-b-c$ and $u(b) \leq (a)$ then $u(c) \leq u(b)$, and if $u(b) <
u(a)$ then $u(c) < u(b)$.

If $G_v$ contains 2 points of branching $x$ and $y$, $x-z_1-\cdots -z_t-y$
($t \geq 0$) and $z_1,\ldots ,z_t$ are not points of branching then $u(x)
\geq u(z_1) \geq \ldots u(z_t) \geq u(y)$ and $u(y) \geq u(z_1) \geq \ldots
u(z_t) \geq u(x)$. Thus, 2) implies $x$ and $y$ are points of weak
branching and $G = \widetilde{D}_n$.

If there is only one point of branching $z$ in $G$ then it is easy to see
that only the cases $\widetilde{E}_6$, $\widetilde{E}_7$,
$\widetilde{E}_8$ and listed above $u_G$ are possible, and the case
of missing of the points of branching is impossible. $\Box$

\medskip

Note that the vector $u_G$ for extended Dynkin graph $G$ is its
(unique up to common multiplier) imaginary root. Clearly,
$u$ is regular.

It is well-known that if $G$ is a Dynkin graph ($A_n$, $D_n$, $E_6$, $E_7$,
$E_8$) then all vectors from $V_G$ are singular \cite{BerGelPon}. Converse
is also true: for all other woods there exist also regular vectors.
Indeed, if $G$~--- extended Dynkin graph then $u_G$ is regular. From
the root theory of the arbitrary graph \cite{Kac2} it is follows
that any imaginary root is regular (we will obtain it also as a corollary
from the proposition~\ref{p_sin_root}). However, also real roots can be
regular: for instance, if $G = \widetilde{E}_6$, $v \in V_G^+$, $v_i = 1$
for $i = \overline{1,7}$.

Locally-scalar representation $\pi$ of graph a $G$ is {\it singular} if
$\pi$ is indecomposable and finite-dimensional, and $d(\pi)$ is a
singular vector; {\it regular} if $\pi$ is indecomposable,
finite-dimensional and not singular.

\begin{hypothesis}\label{h2}

If $G$ is an extended Dynkin graph, $\pi$~--- its regular
representation then $d(\pi) \leq u_G$.

\end{hypothesis}

In \cite{KrugRoit} for (arbitrary) graph $G$ were considered: a category
${\rm Rep}\,(G,\mathcal{H})$ of representations of $G$ in the category
of Hilbert spaces, a category ${\rm Rep}\,(G)$ of locally-scalar
representations, its (full) subcategory ${\rm Rep}\,(G,d,f)$ of
representations with fixed dimension $d$ and character $f$, and,
finally, ${\rm Rep}\,(G,\coprod)$~--- a union of the categories ${\rm
Rep}\,(G,d,f)$ satisfying the condition $d_i+f_i > 0$, $i \in
\overline{1,n}$.

Note that all listed categories, except ${\rm Rep}(G,\mathcal{H})$,
are not additive.

If vector $f$ is such that $f_i > 0$ when $i = \overline{p+1,n}$ then
in \cite{KrugRoit} it was constructed a functor $\iiSC{F_{df}}:
{\rm Rep}\,(G,d,f) \to {\rm Rep}\,(G,d^+,\iiSC{f}_d)$ which is equivalence
($d^+ = \iiSC{c}(d)$). At that $d^+_i = d_i$ when $g_i \in
\iiSB{G}$, $d^+_i = \sigma_i(d)$ when $g_i \in \iiSC{G}$; $(\iiSC{f}_d)_i =
f_i$ when $g_i \in \iiSC{G}$, and if $g_i \in \iiSB{G}$ then
$(\iiSC{f_d})_i = \sigma_i(f)$ when $d_i \neq 0$ and $(\iiSC{f_d})_i = f_i$
when $d_i = 0$. Analogously functors $\iiSB{F_{df}} : {\rm
Rep}\,(G,d,f) \to {\rm Rep}\,(G,d^-,\iiSB{f}_d)$ are constructed ($d^- =
\iiSB{c}(d)$).

For $d \in V_G^+$ define $G^d = \{g \in G_v | d(g) > 0,\}$.
We will need slightly another construction of functors $\iiSC{\Phi}_{df}:
{\rm Rep}\,(G,d,f) \to {\rm Rep}\,(G,d^+,\iiSC{f'}_d)$ on condition that
$f(g_i) > 0$ for $g_i \in (G^d \bigcup M(G^d)) \bigcap \iiSC{G}$. These
functors, as functors $\iiSC{F}_{df}$, results from functors
$\iiSC{F}_{X,\delta}$ by the selection of other $\delta$ (corresponding to
$\iiSC{f'}_d$). Here $(\iiSC{f'}_d)_i = f_i$ when $g_i \in \iiSC{G}$, and
if $g_i \in \iiSB{G}$ then when $d_i = 0$, $g_i \in M(G^d)$ or $d^+_i =
0$, $g_i \in M(G^{d^+})$ we define $(\iiSC{f'}_d)_i = f_i$, and in other
cases $(\iiSC{f'}_d)_i = \sigma_i(f)$.

Functoriality of $\iiSC{\Phi}_{df}$ folows from the reasonings
analogous to mentioned in \cite{KrugRoit} during the proof of
functoriality of $\iiSC{F_{df}}$ from ${\rm Rep}\,(G,d,f)$ to ${\rm
Rep}\,(G,d^+,\iiSC{f}_d)$).

Analogously, of course, functor $\iiSB{\Phi}_{df}$ is constructed.

Let $S' = \{(d,f) \in Z_G \times V_G \,|\, f(g) > 0$ when $g \in
M(G^d)\}$,
$\iiSC{S'} = \{(d,f) \in S' \,|\, f(g) > 0$ when $g \in
G^d \bigcap \iiSC{G}\}$,
$\iiSB{S'} = \{(d,f) \in S' \,|\, f(g) > 0$ when $g \in
G^d \bigcap \iiSB{G}\}$.

Let us construct a category ${\rm Rep}(G,\coprod') =
\coprod\limits_{(d,f) \in S'}{\rm Rep}(G,d,f)$. So, objects of
${\rm Rep}(G,\coprod')$ are pairs $(\pi,f)$ where $\pi$ is a
locally-scalar representation, $f$ is its character (which is completely
determined by the representation $\pi$ only if $\pi$ is faithful; dimension
$d$, of course, is always determined by the representation $\pi$).
Morphisms between objects from ${\rm Ob\, Rep}(G,d_1,f_1)$ and ${\rm Ob\,
Rep}(G,d_2,f_2)$ match with morphisms in ${\rm Ob\, Rep}(G,d_1,f_1)$ when
$(d_1,f_1) = (d_2,f_2)$ and are missing on the pairs of object such that
$(d_1,f_1) \neq (d_2,f_2)$.

A category ${\rm Rep}_\circ(G,\coprod') \subset {\rm Rep}(G,\coprod')$ is a
full subcategory of $f$-representations, for which $(d,f)$
such that $f(g) > 0$ when $g \in (M(G^d) \bigcup G^d)
\bigcap \iiSB{G}$ (analogously a category ${\rm Rep}_\bullet(G,\coprod')$
is defined).

It is easy to check that if $(d,f) \in \iiSC{S'}$ then $(d^+,\iiSC{f'})
\in \iiSC{S'}$ and if $(d,f) \in \iiSB{S'}$ then $(d^+,\iiSB{f'}) \in
\iiSB{S'}$. Thus, a functor $\iiSC{\Phi}
: {\rm Rep}_\circ(G,\coprod') \to {\rm Rep}_\circ(G,\coprod')$~---
a union of functors $\iiSC{\Phi}_{df}$, and a functor
$\iiSB{\Phi}
: {\rm Rep}_\bullet(G,\coprod') \to {\rm Rep}_\bullet(G,\coprod')$~---
a union of functors $\iiSB{\Phi}_{df}$ are determined (analogously
$\iiSC{F}$ and $\iiSB{F}$ in \cite{KrugRoit}).

$\iiSC{\Phi^2} \cong {\rm Id}$ and
$\iiSB{\Phi^2} \cong {\rm Id}$. Therefore, $\iiSC{\Phi}$ ($\iiSB{\Phi}$) is
an equivalence in ${\rm Rep}_\circ\,(G,\coprod')$ (${\rm
Rep}_\bullet\,(G,\coprod')$).

Let us introduce denotations: $\iiSB{c}= \sigma_{g_p}\cdots \sigma_{g_1}$,
$\iiSC{c} = \sigma_{g_n}\cdots \sigma_{g_{p+1}}$. Let $c_t =
\underbrace{\cdots \iiSB{c}\iiSC{c}\iiSB{c}}_{t}$ when $t > 0$,
$c_t = \underbrace{\cdots \iiSC{c}\iiSB{c}\iiSC{c}}_{t}$ when $t < 0$ and
$c_0 = {\rm id}$.

Operators $c_i$ generate in the group of invertible operators in $V_G$
a subgroup isomorphous to dihedral group; $c_r c_s = c_t$, where $t = (-1)^s
r + s$.

Let $\iiRGUn{1} = \iiRGUC \bigcap \iiRGUB$, functors $\iiSC{\Phi}$,
$\iiSB{\Phi}$ are defined on $\iiRGUn{1}$ with values in $\iiRGU$. Let us
construct in $\iiRGU$ a sequence of full subcategories

$$\iiRGU \supset \iiRGUn{1} \supset \ldots \supset \iiRGUn{k} \supset
\ldots$$

\noindent such that ${\rm Ob}\,\iiRGUn{i+1} = \{X \in {\rm
Ob}\,\iiRGUn{i}\,|\, \iiSC{\Phi} \in \iiRGUn{i}, \iiSB{\Phi} \in
\iiRGUn{i}\}$. Thus, on $\iiRGUn{k}$ functors
$\iiSC{\Phi}$, $\iiSB{\Phi}$ with values in $\iiRGUn{k-1}$ are defined.
Hence for each $k \in \mathbb{N}$ functors

$$\Phi_k : \iiRGUn{k} \to \iiRGU,\;\, \Phi_k =
\underbrace{\cdots \iiSB{\Phi}\iiSC{\Phi}\iiSB{\Phi}}_{k},$$

$$\Phi_{-k} : \iiRGUn{k} \to \iiRGU,\;\, \Phi_{-k} =
\underbrace{\cdots \iiSC{\Phi}\iiSB{\Phi}\iiSC{\Phi}}_{k}$$

\noindent are defined.

Formulas listed above imply that if $\Phi_t(\pi,f) = (\pi_t,f_t)$ then
$d(\pi_t) = c_t(d(\pi))$.

Let $G$~--- extended Dynkin graph, $u_G = (u_k)$. Let us construct a
linear form

$$L_G(x) = \sum\limits_{g_i \in \iiSB{G}}u_ix_i - \sum\limits_{g_j \in
\iiSC{G}}u_jx_j,\;\, x \in V_G.$$

Following statement takes place:

\begin{lemma}\label{p_lin_form}

For any $x \in V_G$ $L_G(c_1(x)) = L_G(c_{-1}(x)) = -L_G(x).$

\end{lemma}

{\bf Proof.}

Denote as $S(x_i) = (\sigma_{g_i}(x))_i + x_i$.

$L_G(c_1(x)) + L_G(x) =
\sum\limits_{g_i \in \iiSB{G}}u_i(\sigma_{g_i}(x))_i -
\sum\limits_{g_i \in \iiSC{G}}u_ix_i +
\sum\limits_{g_i \in \iiSB{G}}u_ix_i -
\sum\limits_{g_i \in \iiSC{G}}u_ix_i =
\sum\limits_{g_i \in \iiSB{G}}u_i((\sigma_{g_i}(x))_i+x_i) -
2\sum\limits_{g_i \in \iiSC{G}}u_ix_i =
\sum\limits_{g_i \in \iiSB{G}}u_iS(x_i) -
2\sum\limits_{g_i \in \iiSC{G}}u_ix_i =
\sum\limits_{g_i \in \iiSC{G}}x_i(S(u_i)-2u_i) =
\sum\limits_{g_i \in \iiSC{G}}x_i((\sigma_{g_i}(u_i))_i-u_i) = 0.$

A case $L_G(c_{-1}(x))$ is examined analogously. $\Box$

\medskip

Let $u = \iiSB{u} + \iiSC{u}$ where $(\iiSB{u})_i = u_i$ when $i \leq p$ and
$(\iiSC{u})_i = u_i$ when $i > p$. Let then $u^{|k|}_\circ =
c_k(\iiSC{u})$ when $k \geq 0$ and $u^{|k|}_\bullet = c_k(\iiSB{u})$ when $k
\leq 0$ and, finally, $\widetilde{u}^k = \alpha u^k$ where $\min\limits_{i =
\overline{1,n}}(\widetilde{u}^k)_i = 1$.

Character of the
representation $\pi$ is said to be {\it odd} (resp. {\it even})
{\it standard} if it equals $\widetilde{u}^k_\bullet$ (resp.
$\widetilde{u}^k_\circ$). Characters $\widetilde{u}^0_\circ$ and
$\widetilde{u}^0_\bullet$ are said to be {\it primary standard}.

Locally-scalar $f$-representation of an extended Dynkin graph $G$ is said to
be {\it standard} if $f$ is a standard character.

In \cite{NazRoit} (increasing) function $\rho$: $\rho(n) =
1+\frac{n-1}{n+1}$, $n \in \mathbb{N}_0$ is defined. (It is naturally to
consider $\rho(\infty) = 2$). We will give an explicit formula for
the standard characters in terms of $\rho$ and will show for each
singular root exactly one standard representation corresponds.

\begin{proposition}\label{p_st_rep}

Let $G$~--- extended Dynkin graph. Let $\widetilde{u}^0_\bullet =
(w_1, \ldots ,w_p, 0,\ldots ,0)$ and $\widetilde{u}^0_\circ =
(0,\ldots ,0,w_{p+1}, \ldots ,w_n)$. Then

$\widetilde{u}^{2m}_\bullet =
(w_1, \ldots w_p, \rho(2m)w_{p+1}, \ldots \rho(2m)w_n),$

$\widetilde{u}^{2m+1}_\bullet =
(w_1, \ldots w_p, (4-\rho(2m))w_{p+1}, \ldots (4-\rho(2m))w_n),$

$\widetilde{u}^{2m}_\circ =
(w_1, \ldots w_p, (4-\rho(2m-1))w_{p+1}, \ldots (4-\rho(2m-1))w_n),$

$\widetilde{u}^{2m+1}_\circ =
(w_1, \ldots w_p, \rho(2m+1)w_{p+1}, \ldots \rho(2m+1)w_n),$

$m \in \mathbb{N}_0.$

\end{proposition}

{\bf Proof} is carried out separately for each extended Dynkin graph.
We will illustrate the proof by the example of graph $\widetilde{E}_7$ for
character $\widetilde{u}^0_\bullet$.

(In the picture numbers mean order numbers of vertices.)

\medskip

\centerline{\includegraphics{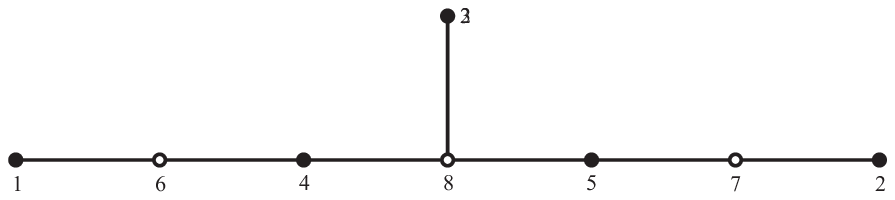}}

\medskip

In this case $\widetilde{u}^0_\bullet = (1,1,2,3,3,0,0,0)$,
$\widetilde{u}^0_\circ = (0,0,0,0,0,1,1,2)$.

Let us carry out an induction for $l = 2m$.

The base of induction is evident.

Let $c_{2m}(\widetilde{u}^0_\bullet) =
(1,1,2,3,3,\rho(2m),\rho(2m),2\rho(2m))$. Hence
$c_{2m+2}(\widetilde{u}^0_\bullet) =
(3-\rho(2m),3-\rho(2m),6-2\rho(2m),9-3\rho(2m),9-3\rho(2m),4-\rho(2m),4-
\rho(2m),8-2\rho(2m))$. Dividing all coordinates by $3-\rho(2m)$ we will
obtain a vector $(1,1,2,3,3,\rho(2m+2),\rho(2m+2),2\rho(2m+2))$.

Now $\widetilde{u}^{2m+1}_\bullet = \iiSC{c}(\widetilde{u}^{2m}_\bullet)
= (1,1,2,3,3,(4-\rho(2m)),(4-\rho(2m)),(8-2\rho(2m)))$.

\medskip

Object $(\pi,f) \in {\rm Rep}\,(G,\coprod')$ is {\it singular} if
$\pi$ is singular. Simplest object in the category ${\rm
Rep}\,(G,\coprod')$ is a pair $(\Pi_g,\bar{f})$ ($d(\Pi_g) = \bar{g}$,
$g \in G_v$, then $\bar{f}(g)$ = 0 and $\bar{f}(x) > 0$ when $x \in M_g$).

\begin{proposition}\label{p_ob_fr_sim}

Each singular object of the category ${\rm Rep}\,(G,\coprod')$
can be obtained in the form of $\Phi_m(\Pi_g,\bar{f})$ where
$(\Pi_g,\bar{f})$~--- simplest object ($m \geq 0$ when $g \in \iiSC{G}$ and
$m \leq 0$ when $g \in \iiSB{G}$). At that each faithful singular
representation $G$ corresponds (up to equivalence) to
unique singular object ${\rm Rep}\,(G,\coprod')$.

\end{proposition}

{\bf Proof} is carried out by the induction on $|t|$ where
$t(d(\pi))$ is minimal in absolute value number, for which
$c_t(d(\pi)) \not\in V_G^+$.

Let $|t| = 1$. For definiteness $\iiSC{c}(d(\pi)) \not\in V_G^+$.
Hence $\exists g \in G^\pi \bigcap \iiSC{G}$, $f(g) = 0$ (in the contrary
case we can apply to $\pi$ the functor $\iiSC{F}_{X,\delta}$ where $X =
G^\pi$, $\delta$ is arbitrary, and will obtain $d(\iiSC{F}_{X,\delta}(\pi))
= \iiSC{c}(d(\pi)) \in V_G^+$). Then indecomposability of $\pi$ and
lemma~3.5 \cite{KrugRoit} imply $G^\pi = \{g\}$, i.~e. $\pi =
\Pi_g$, $(\pi,f)$ is a simplest object.

Let the statement is true with $|t| < k$. Let $(\pi,f)$~--- singular
object ${\rm Rep}\,(G,\coprod')$, $c_k(d(\pi)) \not\in V_G^+$; since
$\pi$ is not a simplest representation then $f$ is positive on $G^\pi$.
Hence either $c_{k-1}(d(\iiSC{F}_{X,\delta}(\pi))) \not\in V_G^+$ or
$c_{k-1}(d(\iiSB{F}_{X,\delta}(\pi))) \not\in V_G^+$ and we cane make use
of the inductive presumption. $\Box$

\begin{corollary}\label{p_sin_root}

Let $x$ is a singular root of extended Dynkin graph $G$. Then $x$ is a real
root in $G$.

\end{corollary}

\begin{proposition}\label{p_form}

Let $G$~--- extended Dynkin graph. In order to root $x \in V_G$
to be singular it is necessary and sufficiently that $L_G(x) \neq 0$.

\end{proposition}

{\bf Proof.}

{\bf Necessity.} If $x$ is singular then by the
proposition~\ref{p_ob_fr_sim} $x = c_t(\bar{g})$ where $\bar{g}$ is a
simple root in $G$. Since $L(\bar{g}) \neq 0$ then, using
the lemma~\ref{p_lin_form}, we will obtain the required.

{\bf Sufficiency.} Consider along with $L_G(x)$ a linear form
$L^+_G(x) = \sum\limits_{g_i \in \iiSB{G}}u_ix_i +
\sum\limits_{g_j \in \iiSC{G}}u_jx_j$. Reckon for definiteness
$L_G(x) > 0$. Since $L_G(c_1(x)) = -L_G(x)$ (lemma~\ref{p_lin_form}) then
$\sum\limits_{g_j \in \iiSB{G}}u_jx_j - \sum\limits_{g_j \in
\iiSB{G}}u_j(\sigma(x_j))_j > 0$, and, consequently, $L^+_G(c_1(x)) <
L^+_G(x)$. Analogously $L^+_G(c_2(x)) < L^+_G(c_1(x))$ etc. Thus, for the
finite number of steps $t$ we will obtain $L^+_G(x) < 0$, which means a
singularity of $x$. $\Box$

\begin{proposition}\label{p_uniq}

If $G$~--- extended Dynkin graph then the presentation of each
singular object of the category ${\rm Rep}\,(G,\coprod')$ in the form of
$\Phi_m(\Pi_g,\bar{f})$ is unique. ($(\Pi_g,\bar{f})$ is a simplest
object, $m \geq 0$ when $g \in \iiSC{G}$ and $m \leq 0$ when $g \in
\iiSB{G}$).

\end{proposition}

It follows from propositions~\ref{p_ob_fr_sim} and \ref{p_form}. $\Box$

\begin{proposition}\label{p_ending}

If $G$~--- extended Dynkin graph, $d$~--- singular root in $G$
then there exists a unique standard representation $\pi$ with the
dimension $d$.

\end{proposition}

{\bf Proof.}

Proposition~\ref{p_ob_fr_sim} imply $x = c_t(\bar{g})$. Let for
definiteness $t \geq 0$. Hence $F_t(\Pi_g,\widetilde{u}^0_\bullet) =
(\pi,\alpha\widetilde{u}^t_\bullet)$. Let $\pi' =
\frac{1}{\sqrt{\alpha}}\pi$ (i.~e. we multiply by $\frac{1}{\sqrt{\alpha}}$
each operator in the representation $\pi$). Then
$\widetilde{u}^t_\bullet$ is a character of the representation $\pi'$.

Uniqueness $\pi$ follows from the proposition~\ref{p_uniq}. $\Box$

%**************************************************************************

\section{Generalization of function $\rho$.}

Let $\widetilde{V}$~--- set of infinite nonincreasing
sequences of the integer nonnegative numbers
$v = (v_1,v_2,\ldots,v_i,\ldots)$ such that $v_s = 0$ from the certain
$s$.  Define on $\widetilde{V}$ a partial order: $v \leq w$,
if $v_i \leq w_i$ for all $i \in \mathbb{N}$.

Function $\rho$ can be defined on set $\widetilde{V}$:
for $v \in \widetilde{V}$ $\rho(v) = \sum\limits_{i=1}^m \rho(v_i)$,
$\rho(0) = 0$. Hence if $v \leq w$ then $\rho(v) \leq \rho(w)$.

{\it Width} $\omega(v)$ of a vector $v \in \widetilde{V}$ is a number of
its nonzero components.

Let us introduce two following lists of vectors from $\widetilde{V}$:

$K = \{(1,1,1,1),\; (2,2,2)\; (3,3,1)\; (5,2,1)\}$;

$\widehat{K} = \{(1,1,1,1,1),\; (2,1,1,1)\; (3,2,2)\; (4,3,1)\;
(6,2,1)\}$.

(Here and further, writing vectors from $\widetilde{V}$, we will write
only their nonzero components.)

In the applications \cite{NazRoit} one often encounter with
conditions

\begin{equation}\label{f21}
\rho(v) = 4,
\end{equation}

\begin{equation}\label{f22}
\rho(v) > 4.
\end{equation}

\begin{proposition}\label{p_list_K}
\cite{NazRoit}

All solutions of the equation~(\ref{f21}) are exhausted by the list $K$; all
minimal solutions of the inequality~(\ref{f22}) are exhausted by the list
$\widehat{K}$.

\end{proposition}

\begin{proposition}\label{p_wv}

If $\rho(w) > 4$ for $w \in \widetilde{V}$ then $\exists v \in
\widetilde{V}\, |\, v < w \mbox{ and }\rho(v) = 4$.

\end{proposition}

{\bf Proof} follows from that fact that for any $\widehat{v} \in
\widehat{K}$ ($\widehat{v} \leq w$) there exists such $v \in K$ that $v <
\widehat{v}$. $\Box$

\medskip

Let us find out, what in propositions \ref{p_list_K} and \ref{p_wv}
is elementary and what is connected with the specificity of function $\rho$.

Let $\varphi : \mathbb{N} \to \mathbb{R}^+$ is an arbitrary function.
We will say that $\varphi$ is {\it convex} if $\varphi(n+1) -
\varphi(n) < \varphi(m+1) - \varphi(m)$ when $n > m$ and {\it
normalized} if $\varphi(1) = 1$. ($\rho$ possess these properties.)
Let $\varphi(z_1,\ldots,z_t) = \sum\limits_{i=1}^{t}\varphi(z_i)$,
$K(\varphi,n) = \{v \in \widetilde{V} \,|\, \varphi(v) = n\}$ and
$N(\varphi,k)$ to consist of the minimal solutions of inequality
$\varphi(v) > k$. For each $x \in \widetilde{V}$ correspond $\widehat{x} \in
\widetilde{V}$ where $\widehat{x}_1 = x_1+ 1$ and $\widehat{x}_i = x_i$
when $i \neq 1$. Let ${\langle s \rangle}^r =
\underbrace{s,\ldots,s}_{r}$. Let $X \subset \widetilde{V}$, $|X| <
\infty$, $\omega(X) = \max\limits_{x \in X}\omega(x)$ and
$\widehat{X} = \{\widehat{x}\,|\, x \in X\} \bigcup
\{\iiCor{1}{\omega(X)+1}\}$.

Following proposition is proved elementary.

\begin{proposition}

Let $n \in \mathbb{N}$. If $\varphi$ is increasing then $K(\varphi,n)$ and
$N(\varphi,n)$ are finite. If, besides, $\varphi$ is convex and
normalized then $\widehat{K}(\varphi,n) \subset N(\varphi,n)$, and at that
each vector from $\widehat{K}(\varphi,n)$ larger than precisely one
vector from $K(\varphi,n)$. If $\widehat{K}(\varphi,n) = N(\varphi,n)$ and
$m < n$ then $\widehat{K}(\varphi,m) = N(\varphi,m)$.

\end{proposition}

Increasing, convex and normalized function $\varphi$ is said to be {\it
$n$-separating} ($n \in \mathbb{N}$) if $N(\varphi,n) =
\widehat{K}(\varphi,n)$. It is equal to the following: if $\varphi(w) > n$,
$w \in \widetilde{V}$ then there exists $v < w$ such that $\varphi(v) = n$.

Thus, $\rho$ is 4-separating.

Let $v \in \widetilde{V}$, $t \in \mathbb{N}$. Let $v^{\langle t
\rangle} = (v^{\langle t \rangle})_i$ where $(v^{\langle t \rangle})_i =
v_1$ when $i \leq t$, and $(v^{\langle t \rangle})_i = v_{i-t}$ when $i >
t$. For $X \subset \widetilde{V}$ let $X^{\langle t \rangle} =
\{x^{\langle t \rangle}\}$ when $x \in X$.

Let us show that when $n > 4$ $\rho$ is not
$n$-separating and for $n = 4+i$ we construct $n$-separating functions
$\rho_i$, such that $K(\rho_i,4+i) = K(\rho,4)^{\langle i \rangle}$ which
(see\S~3) also connected with locally-scalar representations of quivers.

Consider an inequality

\begin{equation}\label{f23}
\rho(v) \geq 5.
\end{equation}

Distinguish several cases.

1. $\omega(v) > 5$. Minimal vector satisfying this condition
is vector $w = (1,1,1,1,1,1), \rho(w) > 5$. Since when $v < w$
$\rho(v) \leq 5$ then $w \in N(\rho,5)$.

2. $\omega(v) = 5$. Under this condition $\rho(v) \geq 5$ and $\rho(v) = 5$
only if $v = (1,1,1,1,1)$ ($v$~--- the lowest vector in this case), i.~e.
$(1,1,1,1,1) \in K(\rho,5)$. Then, obviously, $(2,1,1,1,1) \in
N(\rho,5)$.

3. $\omega(v) = 4$. Let $v = (v_1,v_2,v_3.v_4)$. If $v_4 \geq 2$ then
$\rho(v) \geq 4\rho(2) = 5\frac{1}{3} > 5$. Hence $(2,2,2,2) \in
N(\rho,5)$. If $v_4 = 1$, $v_3 \geq 2$ then $\rho(v) \geq
\rho(1)+3\rho(2) = 5$, in this case we obtain $(2,2,2,1) \in K(\rho,5)$,
$(3,2,2,1) \in N(\rho,5)$. If $v_3=v_4=1$ then the inequality reduces to
$\rho(v_1)+\rho(v_2) \geq 3$. When $v_2 = 2$ we obtain $(5,2,1,1) \in
K(\rho,5)$ and $(6,2,1,1) \in N(\rho,5)$. When $v_2 \geq 3$ we obtain
$(3,3,1,1) \in K(\rho,5)$ and $(4,3,1,1) \in N(\rho,5)$.

4. $\omega(v) = 3$. Let $v = (v_1,v_2,v_3)$.
Using reasonings, analogous to previous, we obtain:

a) when $v_3 \geq 5$ $(5,5,5) \in K(\rho,5)$ and $(6,5,5) \in N(\rho,5)$;

b) when $v_3 = 4$ $(9,4,4) \in K(\rho,5)$ and $\{(10,4,4), (7,5,4)\} \subset
N(\rho,5)$;

c) when $v_3 = 3$ $\{(19,4,3),(11,5,3)\} \subset K(\rho,5)$ and $\{(20,4,3),
(12,5,3), (9,6,3)\} \subset N(\rho,5)$;

d) when $v_3 = 2$ $\{(41,6,2), (23,7,2), (17,8,2), (14,9,2), (11,11,2)
\subset K(\rho,5)\}$; $\hfill$ $\{(42,6,2),
(24,7,2), (18,8,2), (15,9,2), (12,11,2), (13,10,2) \subset N(\rho,5)\}$;

e) when $v_3 = 1$ we have $\rho(v) < \rho(1) + 2\rho(\infty) = 5$.

5. The case $\omega(v) \leq 2$ does not give any solutions, because
$\rho(v) < 2\rho(\infty) = 4 < 5$.

Thus, the following proposition takes place:

\begin{proposition}

$K(\rho,5) = $
$\{(1,1,1,1,1),$ $(2,2,2,1),$ $(3,3,1,1),$ $(5,2,1,1),$  $(41,6,2),$
$(23,7,2),$ $(17,8,2),$ $(14,9,2),$ $(11,11,2),$ $(19,4,3),$ $(11,5,3),$
$(7,7,3),$ $(9,4,4),$ $(5,5,5).\}$. $N(\rho,5) = \widehat{K}(\rho,5)
\bigcup \{(9,6,3), (7,5,4), (2,2,2,2), (13,10,2)\}$

\end{proposition}

Let us fix $t \in \mathbb{N}$ and define a recurrent
sequence $\{u_i\}$ in the following way:

\begin{equation}\label{f25}
u_0 = 0, u_1 = 1, u_{i+2} = tu_{i+1}-u_i.
\end{equation}

Then, let us define a function $\rho_{t-2}(n)$:

\begin{equation}\label{f26}
\rho_{t-2}(n) = 1+\frac{u_{n-1}}{u_n + 1} \mbox{ when } n \in \mathbb{N};\;\,
\rho_{t-2}(0) = 0.
\end{equation}

With $t=2$ from (\ref{f25}) it is follows that $u_n = n$ and $\rho_0(n) =
1+\frac{n-1}{n+1} = \rho(n)$, i.~e. $\rho_0(n) = \rho(n)$.

Let then $k = t-2$, $k \in \mathbb{N}_0$.

Let us find a general formula for $n$-th member of sequence
$\{u_n\}$. For that we solve a characteristic equation

$$\lambda^2-(k+2)\lambda+1=0.$$

Its roots are $\lambda = \frac{k+2+\sqrt{k^2+4k}}{2}$, $\overline\lambda =
\frac{k+2-\sqrt{k^2+4k}}{2}$.

Then formula of the $n$-th member (when $t \neq 2$) has the form

$$u_n = C_1\lambda^n + C_2\overline{\lambda}^n.$$

From the entry conditions we obtain $C_1 = \frac{1}{\sqrt{k^2+4k}},\; C_2
= -\frac{1}{\sqrt{k^2+4k}}.$

Finally

$$u_n = \frac{1}{\sqrt{k^2+4k}} (\lambda^n - \overline{\lambda}^n).$$

After trivial transformations from the formula~(\ref{f26}) we obtain

\begin{equation}\label{f27}
\rho_k(n) = 1 + \frac{\lambda^n - \lambda}{\lambda^{n+1} - 1},\; k \neq 0.
\end{equation}

Let us calculate

\begin{equation}\label{f28}
\rho_k(\infty) = \lim\limits_{n \to \infty}\rho_k(n) = 1 +
\overline{\lambda} = 1 + \frac{2}{k+2+\sqrt{k^2+4k}}.
\end{equation}

Let us also show that the function $\rho_k(n)$ is increasing on
$\mathbb{N}$. Using formula~(\ref{f27}), the inequality $\rho_k(n) <
\rho_k(n+1)$ reduces to the equivalent one $\lambda^3 + 1 > \lambda^2 +
\lambda$ which is true due to $\lambda > 1$.

Define with $k > 0$ a function $\rho_k(n)$ on $\widetilde{V}$ analogously
to the function $\rho$: for $v \in \widetilde{V}$ let $\rho_k(v) =
\rho_k(v_1,\ldots,v_s) = \sum\limits_{i=1}^s \rho_k(v_i)$. Since
$\rho_k$ is increasing on $\mathbb{N}_0$ it is increasing also on
$\widetilde{V}$.

Consider an inequality

\begin{equation}\label{f29}
\rho_k(v) \geq k+4
\end{equation}

Let us show that the results will be analogous to results obtained
for function $\rho$.

\begin{lemma}\label{l_for}

For any $k,n \in \mathbb{N}_0$ the following inequality is correct:

\begin{equation}\label{f211}
\rho_k(n) + k\rho_k(\infty) < (k+1)\rho_k(n+1)
\end{equation}

\end{lemma}

{\bf Proof.}

When $k = 0$ the inequality reduces to the evident one $\rho(n) <
\rho(n+1)$. Let $k > 0$. Using formulas (\ref{f27}) and (\ref{f28}), after
transformations we obtain that the inequality (\ref{f211}) is
equivalent to inequality

$$(\lambda^2 - 1)(\lambda^{n+2} - (k+1)\lambda^{n+1} + k) > 0,$$

\noindent which is true since
$\lambda = \frac{k+2+\sqrt{k^2+4k}}{2} > k+1$ when $k > 0$. $\Box$

\medskip

\begin{proposition}\label{p_list_Kt}

Function $\rho_k$ is $(k+4)$-separating, $K(\rho_k,4+k) = K^{\langle k
\rangle}$ $(K = K(\rho,4))$.

$K(\rho_k,4+k) = \{(\iiCor{1}{k+4});\; (\iiCor{2}{k+3});\;
(\iiCor{3}{k+2},1);\; (\iiCor{5}{k+1},2,1)\}$;

\end{proposition}

{\bf Proof.}

On the ground of formula~(\ref{f26}) we will find sequentially $\rho_k(1) =
1$, $\rho_k(2) = 1+\frac{1}{k+3}$, $\rho_k(3) = 1+\frac{1}{k+2}$,
$\rho_k(4) = 1+\frac{k+3}{k^2+5k+5}$, $\rho_k(5) = 1+\frac{k+2}{k^2+4k+3}$.

Consider inequality (\ref{f29}).

1. Let $\omega(v) > k+4$. Then $\rho_k(v) \geq
\rho_k(\iiCor{1}{k+5})$. Minimal vector among vectors with this width
is vector $w = (\iiCor{1}{k+5})$. It is clear that for any $v <
w$ (\ref{f29}) is correct, therefore $w \in N(\rho_k,k+4)$.

2. $\omega(v) = k+4$. $\rho_k(v) \geq \rho_k(\iiCor{1}{k+4}) = k+4$,
and at that the equality here will be only if $v = (\iiCor{1}{k+4})$.
Therefore, $(\iiCor{1}{k+4}) \in K(\rho_k,k+4)$. Minimal vector in this
case is $w = (2,\iiCor{1}{k+3})$. The same way as in the previous case,
for any $v < w$ will be $\rho_k(v) \leq k+4$, hence $w \in
N(\rho_k,k+4)$.

3. $\omega(v) = k+3$. Let $v = (v_1,\ldots,v_{k+3})$.

a) Assume $v_{k+3} \geq 2$. Hence $\rho_k(v) \geq
(k+3)\rho_k(2) \geq k+4$. Then the equality will be only if $v =
(\iiCor{2}{t+1}) \in K(\rho_k,k+4)$, minimal vector is vector
$w = (3,\iiCor{2}{t}) \in N(\rho_k,k+4)$.

b) Let $v_{k+3} = 1$ and $v_{k+2} \geq 2$. If $v_{k+2} \geq 3$ then
$\rho_k(v) \geq \rho_k(1) + (k+2)\rho_k(3) = k+4$. So:
$(\iiCor{3}{k+2},1 \in K(\rho_k,k+4))$, $(4,\iiCor{3}{k+1},1 \in
N(\rho_k,k+4))$. Now consider $v_{k+2}$ = 2. The initial inequality
reduces to the following: $\rho_k(v_1)+\ldots +\rho_k(v_{k+1}) +
\frac{1}{k+3} \geq {k+2}$. The equality is correct here if all $v_j = 5$,
$j = \overline{1,k+1}$. Let us show that for all other vectors $v =
(v_1,\ldots,v_{k+1},2,1)$ in the case $v_{k+1} < 5$ will be $\rho_k(v) <
\rho_k(\iiCor{5}{k+1},2,1)$. Indeed, by the lemma~\ref{l_for}
$\rho_k(v) < \rho_k(1) + \rho_k(2) + \rho_k(4) +
k\rho_k(\infty) < \rho_k(1) + \rho_k(2) +
(k+1)\rho_k(5) = \rho_k(\iiCor{5}{k+1},2,1)$. And if $v_{k+1} \geq
5$ then vector $(\iiCor{5}{k+1},2,1)$ will be, obviously, minimal.
Thus, in this case we have $(\iiCor{5}{k+1},2,1) \in K(\rho_k,k+4)$ and
$(6,\iiCor{5}{k},2,1) \in N(\rho_k,k+4)$.

c) $v_{k+3} = v_{k+2} = 1$. $\rho_k(v) <
\rho_k(\iiCor{\infty}{k+1},1,1) = 2+(n-1)\rho_k(\infty)$. Hence by
the formula~(\ref{f28}) $\rho_k(v) < k+3 + \frac{k+1}{\lambda} < k+4$.

4. $\omega(v) \leq k+2$. $\rho_k(v) < (k+2)\rho_k(\infty) \leq k+4$
since $\rho_k(\infty) = 1+\frac{1}{\lambda} \leq 1+\frac{2}{k+2}$. In this
case there are no solutions.

We considered all possible cases and, accordingly, we obtained that

$$K(\rho_k,k+4) = \{(\iiCor{1}{k+4}), (\iiCor{2}{k+3}),
(\iiCor{3}{k+2},1), (\iiCor{5}{k+1},2,1)\};$$

$$N(\rho_k,k+4) = \{(\iiCor{1}{k+4}), (2,\iiCor{1}{k+3}),
(3,\iiCor{2}{k+2}), (4,\iiCor{3}{k+1},1), (6,\iiCor{5}{k},2,1)\}.$$

$\Box$

\medskip

Let $\alpha \in \mathbb{R}^+$. Consider functions $\rho_\alpha(v)$,
$v \in \widetilde{V}$ with this assumption, defining them by the same rules
(\ref{f25}) and (\ref{f26}) ($\alpha + 2 = t$). Now let us fix
$v = (v_1,\ldots,v_s)$, $v \in \widetilde{V}$ and consider an equation

\begin{equation}\label{f212}
\rho_\alpha(v) = \alpha + 4
\end{equation}

\noindent where $\alpha$ is unknown.

\begin{proposition}\label{p_irr_t}

Equation~(\ref{f212}) has no rational non-integer solutions.

\end{proposition}

{\bf Proof.}

Formula~(\ref{f25}) by induction evidently implies that $u_n$
has a form $u_n = t^{n-1} + a_{n-2}t^{n-1} + \ldots + a_0$, $a_i \in
\mathbb{Z}$, $i \in \mathbb{N}_0$. Let $v = (v_1,\ldots,v_s)$. Hence
the equation~(\ref{f212}) (taking into account~(\ref{f26})) will have a form

$$\sum\limits_{i=1}^s (1+\frac{u_{v_i-1}}{u_{v_i}+1}) = t+2$$

$$\sum\limits_{i=1}^s (\frac{t^{n-2} + b_{n-3}^{(i)}t^{n-3} + \ldots +
b_0^{(i)}}{t^{n-1} + a_{n-2}^{(i)}t^{n-1} + \ldots + a_0^{(i)} + 1}) =
t+2-s\;\; (a^{(i)}_j,b^{(i)}_j \in \mathbb{Z}).$$

After reduction to the common denominator we obtain an equation $P(t) = 0$
where $P(t)$ is a polynomial with coefficients from $\mathbb{Z}$,
which leading coefficient equals 1. Consequently, $P(t)$ over field
$\mathbb{R}$ has either integer or irrational roots. $\Box$

\begin{remark}

Proposition~\ref{p_irr_t} is not true if assume $v_i \in
\mathbb{N} \bigcup \{\infty\}$. For instance, $\rho_{0,5}(1,1,1,\infty) =
4,5$.

\end{remark}

So, when $\alpha \in \mathbb{Q}^+$, the equality $\rho_\alpha(v) = \alpha +
4$ is possible only with integer $\alpha$ (proposition~\ref{p_irr_t});
for any $\alpha \in \mathbb{N}_0$ there exists a list of vectors
satisfying this equality (proposition~\ref{p_list_Kt}).

%*************************************************************************

\section{Standard representations of star graphs.}

Graph $G$ is said to be a {\it star graph} or a {\it graph of type
$(n_1,\ldots ,n_s)$ ($n_i \in N$, $i = \overline{1,s}$, $s > 1$)}, if $G$
is a wood and $G$ has exactly one vertex of multiplicity $s$ ({\it nodal}),
$s$ vertices of multiplicity 1 ({\it extreme}), and the other vertices have
multiplicity 2, at that the number of edges between nodal vertex and
extreme ones equal $n_1,\ldots ,n_s$. The set of such graphs we denote $S$.

Let now $S_0 \subset S$~--- the set of graph of the type $(n_1,\ldots ,n_s)$
such that $\rho_k(n_1,\ldots ,n_s) = k+4$, $k \in \mathbb{N}_0$. The
constitution of the set $S_0$ is determined by the
proposition~\ref{p_list_Kt}. Note that in $S_0$ extended
Dynkin graphs $\widetilde{D}_n$, $\widetilde{E}_6$, $\widetilde{E}_7$,
$\widetilde{E}_8$ are contained (they correspond to the case $k = 0$).

Let us consider locally-scalar representations of a graph from $S_0$.
Further we will obtain a result analogous to the
proposition~\ref{p_st_rep}. Below we specify for each family of graphs from
$S_0$ even and odd standard characters. (In the pictures numbers mean the
order number of vertices.)

1. $G = (\iiCor{1}{k+4})$.

\medskip

\centerline{\includegraphics{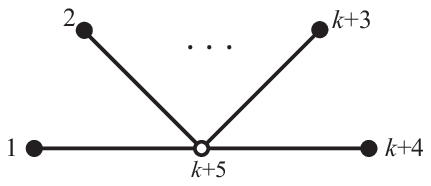}}

\medskip

\noindent $\widetilde{u}^0_\bullet(k) = (\iiCor{1}{k+4},0$);
\noindent $\widetilde{u}^0_\circ(k) = (\iiCor{0}{k+4},1$).

\medskip

2. $G = (\iiCor{2}{k+3})$.

\medskip

\centerline{\includegraphics{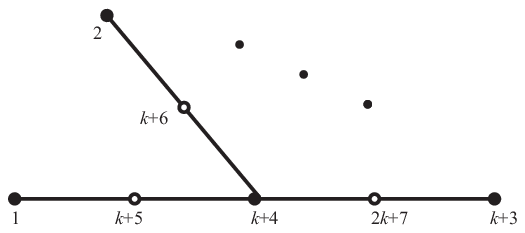}}

\medskip

\noindent $\widetilde{u}^0_\bullet(k) =
(\iiCor{1}{k+3},(k+3),\iiCor{0}{k+3}$);
\noindent $\widetilde{u}^0_\circ(k) = (\iiCor{0}{k+4},\iiCor{1}{k+3})$.

\medskip

3. $G = (\iiCor{3}{k+2},1)$.

\medskip

\centerline{\includegraphics{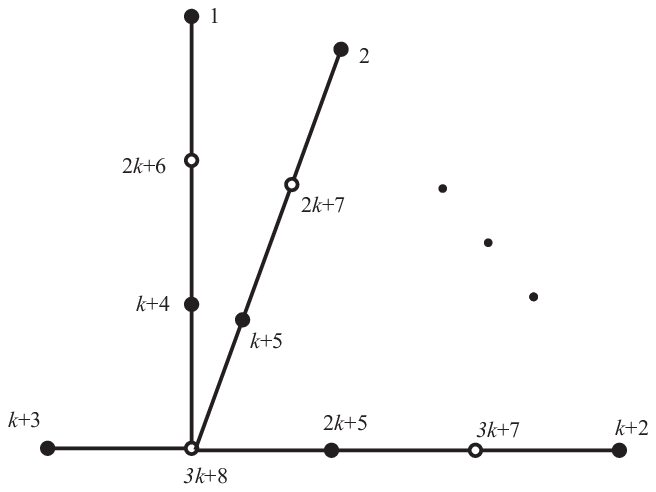}}

\medskip

\noindent $\widetilde{u}^0_\bullet(k) =
(\iiCor{1}{k+2},k+2,\iiCor{k+3}{k+2},\iiCor{0}{k+3})$;
\noindent $\widetilde{u}^0_\circ(k) = (\iiCor{0}{k+2},\iiCor{1}{k+3},k+2)$.

\medskip

4. $G = (\iiCor{5}{k+1},2,1)$.

\medskip

\centerline{\includegraphics{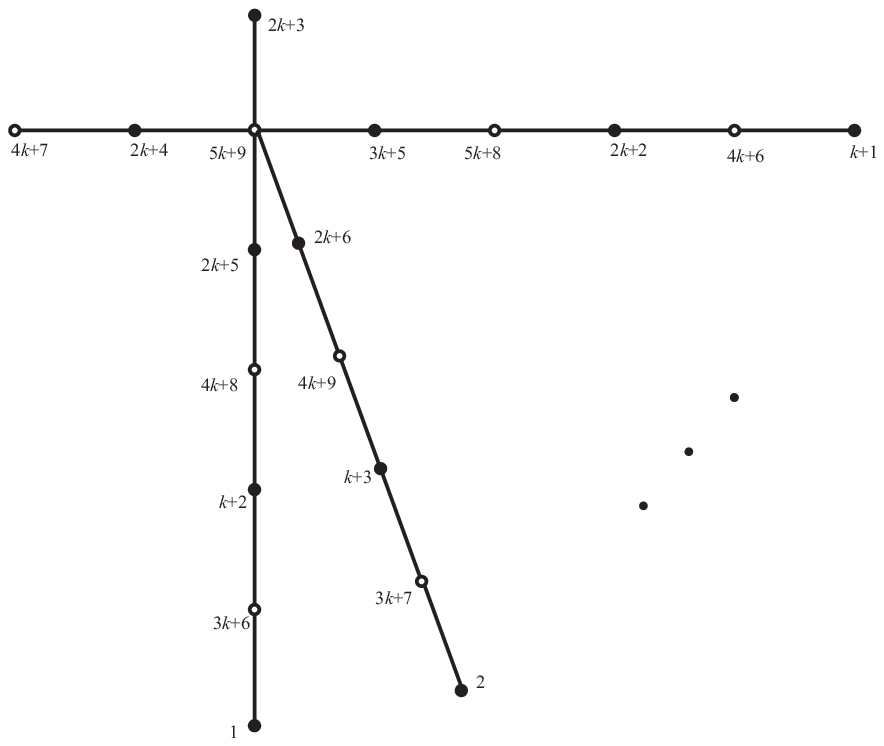}}

\medskip

\noindent $\widetilde{u}^0_\bullet(k) =
(\iiCor{1}{k+1},\iiCor{k+3}{k+1},(k^2+4k+3),(k^2+5k+4),\iiCor{k^2+5k+5}{k+1}
,\iiCor{0}{2k+2})$;
\noindent $\widetilde{u}^0_\circ(k) =
(\iiCor{0}{3k+2},\iiCor{1}{k+1},(k+1),\iiCor{k+2}{k+1},(k^2+4k+3))$.

Hence the following proposition, analogous to the
proposition~\ref{p_st_rep}, takes place.

\begin{proposition}\label{p_st_rep_gen}

Let $G \in S_0$,
$\widetilde{u}^0_\bullet(t) = (w_1, \ldots ,w_p, 0,\ldots ,0)$ and
$\widetilde{u}^0_\circ(t) = (0,\ldots ,0,w_{p+1}, \ldots ,w_n)$. Then:

$\widetilde{u}^{2m}_\bullet(t) =
(w_1, \ldots w_p, \rho_t(2m)w_{p+1}, \ldots \rho_t(2m)w_n),$

$\widetilde{u}^{2m+1}_\bullet(t) =
(w_1, \ldots w_p, (t+2-\rho_t(2m))w_{p+1}, \ldots (t+2-\rho_t(2m))w_n),$

$\widetilde{u}^{2m}_\circ(t) =
(w_1, \ldots w_p, (t+2-\rho_t(2m-1))w_{p+1}, \ldots (t+2-\rho_t(2m-1))w_n),$

$\widetilde{u}^{2m+1}_\circ(t) =
(w_1, \ldots w_p, \rho_t(2m+1)w_{p+1}, \ldots \rho_t(2m+1)w_n),$

$m \in \mathbb{N}_0.$

\end{proposition}

{\bf Proof} is carried out analogous to the proof
of the proposition~\ref{p_st_rep} by the induction on $l = 2m$ separately
for each family of graphs from $S_0$.

Formula for standard characters immediately implies

\begin{proposition}\label{p_stand_inv}

If $f$ is a standard character then $\iiSC{f}$ and $\iiSB{f}$ are
standard character too. If $(\pi,f)$ is a standard object of the category
$\iiRGU$ then $\iiSC{F}(\pi,f)$ is a standard object if $\pi \neq
\Pi_g$, where $g \in \iiSC{G}$, and $\iiSB{F}(\pi,f)$ is a standard
object if $\pi \neq \Pi_g$, where $g \in \iiSB{G}$.

\end{proposition}

\medskip

\begin{proposition}

If $G$ is a star graph and $\rho_k(n_1,\ldots,n_s) = k+4$, $d$~---
singular root in $G$ then there exists an unique singular
$f$-representation $\pi$ with dimension $d$, where $f$ is a standard
character.

\end{proposition}

{\bf Proof} of the existence is analogous to
the proposition~\ref{p_ending}. The uniqueness folows from
the proposition~\ref{p_stand_inv}. $\Box$


\begin{thebibliography}{99}

\bibitem{Gabr} {\it P.~Gabriel.} Unzerlegbare Darstellungen I -- Manuscripta
Math.~6 (1972), 71-103.

\bibitem{BerGelPon} {\it Bernstein~I.N., Gelfand~I.M., Ponomarev~V.A.}
Coxeter functors and Gabriel's theorem. -- IMN, v. XXVIII, pt.~2, p.~19-33
(1973).

\bibitem{Naz2} {\it Nazarova~L.A..} Representations of quivers
of infinite type. -- Izv. AN USSR. Ser. mat. 1973, {\bf 37},
p.~752--791.

\bibitem{DonFre} {\it Donovan~P., Freislich~M.R.} The representation theory
of finite graphs and associated algebras. -- Carleton Math. Lecture Notes,
1973, 5, p.1--119.

\bibitem{Kac} {\it Kac~V.G.} Some remarks on representations of quivers and
infinite root systems. -- Lect. Notes Math., 1980, {\bf 832}, p.~311-332.

\bibitem{GabRoi} {\it P.~Gabriel, A.V.~Roiter} Representations of
Finite-Dimensional Algebras. Berlin etc., Springer, 1997 -- 171 p.

\bibitem{KrugRoit} {\it Kruglyak~S.A., Roiter~A.V.} Locally-scalar
representations of graphs in the category of Hilbert spaces. -- Present
preprint.

\bibitem{OstSam} {\it Ostrovskii~V., Samoilenko~Yu.} Introduction to the
Theory of Representations of Finitely Presented $*$-algebras. I.
Representations by bounded operators. Rev. Math. \& Math. Phys., vol.11,
1--261, Gordon and Breach, 1999.

\bibitem{KruRab} {\it Kruglyak~S.A., Rabanovich~V.I., Samoilenko~Yu.S.} On
sums of projectors -- Functional analysis and its applications, v.~36,
pt.~3, p.~20-35, 2002.

\bibitem{NazRoit} {\it Nazarova~L.A., Roiter~A.V.} Norm of a relation,
separating functions and representations of marked quivers. -- Ukr.
Mat. Jour., 2002, v.~54, N.~6.

\bibitem{Kac2} {\it Kac~V.G.} Infinite root systems, representations of
graphs and invariant theory, II. -- I. Algebra, 1982, {\bf 78}, p.~141-162

\bibitem{GelPon} {\it Gelfand~I.M., Ponomarev~V.A.} Problems of linear
algebra and classification of quadruples of subspaces in a
finite-dimensional vector space. -- Coll. Math. Soc. J. Bolyai, 5. Hilbert
Space Operators, Tihany (Hungary), 1970.

\end{thebibliography}
\end{document}